# A quadratic-residue dichotomy for two partition functions modulo 3

Jiyuan Li
Department of Mathematics
University of California, San Diego
jil537@ucsd.edu

**Abstract.** Let $f_{0,1,4}(n)$ denote the number of partitions of $n$ into parts congruent to 0, 1 or 4 modulo 5 with each part used at most twice, and let $f_{0,2,3}(n)$ be defined analogously for parts congruent to 0, 2 or 3 modulo 5. We prove that for **every** prime $p \equiv 3 \pmod 4$ there are explicit non-negative integers $a(p)$ and $b(p)$, determined by $20\,a(p) \equiv -9$ and $20\,b(p) \equiv -1 \pmod{p^2}$, such that $f_{0,1,4}\big(p^2m + a(p)\big)$ and $f_{0,2,3}\big(p^2m + b(p)\big)$ are congruent modulo 3 to $f_{0,1,4}(m)$ and $f_{0,2,3}(m)$ when $p \equiv \pm 1 \pmod 5$, and to $f_{0,2,3}(m)$ and $f_{0,1,4}(m)$ when $p \equiv \pm 2 \pmod 5$. By quadratic reciprocity the two functions are preserved exactly when 5 is a quadratic residue modulo $p$, and interchanged otherwise. The smallest cases are $f_{0,1,4}(9m) \equiv f_{0,2,3}(m)$ and $f_{0,2,3}(9m+4) \equiv f_{0,1,4}(m)$. The proof reduces each generating function modulo 3 to the square of a Rogers–Ramanujan-type theta function by means of $(1-x)^\ell \equiv 1 - x^\ell \pmod \ell$ and the Jacobi triple product, and then dissects the resulting binary quadratic form using the fact that $-1$ is a quadratic non-residue modulo $p$. As corollaries we obtain, for each such $p$, a congruence with exceptions on a progression of modulus $p$, and an infinite family of self-similarity congruences on progressions of modulus $p^{dn}$ with $d = 2$ or $d = 4$ and constants $9(p^{dn} - 1)/20$ and $(p^{dn} - 1)/20$.

## 1. Introduction

Ramanujan's congruences $p(5n+4) \equiv 0 \pmod 5$, $p(7n+5) \equiv 0 \pmod 7$ and $p(11n+6) \equiv 0 \pmod{11}$ [6] have been the model for a large body of work on arithmetic progressions carrying congruences for restricted partition functions. A particularly active strand concerns the $\ell$-regular partition functions $b_\ell(n)$, whose generating function is the eta-quotient $\left(q^\ell; q^\ell\right)_\infty/(q;q)_\infty$, modulo the small prime 3; see Webb [7] for $\ell = 13$, Furcy and Penniston [3] for general $\ell$, and Keith [5] and Cui and Gu [2] for $\ell = 9$, together with the references in those papers. The method used below is closest in spirit to that of Hou, Sun and Zhang [4], who reduce a generating function modulo 3 to a power of a theta series and then dissect the resulting quadratic form. Keith [5] also introduced the useful notion of a *congruence with exceptions*, a progression on which a partition function is divisible by a prime except on a thin sub-progression; congruences of exactly that shape appear in Section 6.

The two functions studied here differ from the $\ell$-regular family in two respects. The parts are restricted by a congruence condition modulo 5 rather than by a divisibility condition, and their multiplicity is bounded rather than unbounded. The admissible residues are the two Rogers–Ramanujan classes $\pm 1$ and $\pm 2$ modulo 5, each taken together with the multiples of 5; as a consequence the two generating functions reduce modulo 3 to squares of the two Rogers–Ramanujan-type theta series, rather than to powers of eta-quotients, and the two functions are exchanged by some of the dissections below rather than preserved.

**Definition 1.1.** For $n \geq 0$, let $f_{0,1,4}(n)$ denote the number of partitions of $n$ into parts congruent to 0, 1 or 4 modulo 5, each part occurring at most twice, and let $f_{0,2,3}(n)$ denote the number of partitions of $n$ into parts congruent to 0, 2 or 3 modulo 5, each part occurring at most twice. In both cases the empty partition is counted, so $f_{0,1,4}(0) = f_{0,2,3}(0) = 1$.

The first sixteen values, beginning at $n = 0$, are

$$f_{0,1,4}\colon\ 1,\ 1,\ 1,\ 0,\ 1,\ 2,\ 3,\ 2,\ 2,\ 3,\ 6,\ 7,\ 6,\ 5,\ 8,\ 12,\ \dots$$

$$f_{0,2,3}\colon\ 1,\ 0,\ 1,\ 1,\ 1,\ 2,\ 1,\ 3,\ 3,\ 2,\ 6,\ 3,\ 7,\ 7,\ 7,\ 12,\ \dots$$

Our main result is Theorem 5.1: for every prime $p \equiv 3 \pmod 4$ the pair $\left(f_{0,1,4}, f_{0,2,3}\right)$ satisfies a cross or self congruence modulo 3 on an explicit progression of modulus $p^2$, the alternative being decided by whether 5 is a quadratic residue modulo $p$. The two smallest instances are

$$f_{0,1,4}(9m) \equiv f_{0,2,3}(m), \qquad f_{0,2,3}(9m+4) \equiv f_{0,1,4}(m) \pmod 3$$

and

$$f_{0,1,4}(49m+2) \equiv f_{0,2,3}(m), \qquad f_{0,2,3}(49m+22) \equiv f_{0,1,4}(m) \pmod 3.$$

All statements below have been checked numerically for $n \leq 60\,000$.

## 2. Notation

Throughout, $q$ is a formal variable and all series are formal elements of $\mathbb{Z}[[q]]$, so that no question of convergence arises; an infinite product is admissible because each of its coefficients depends on only finitely many factors. We use the $q$-Pochhammer symbol

$$(a;q)_\infty \;=\; \prod_{i\geq 0}\left(1 - aq^i\right).$$

We write $S_1 = \{\, k \geq 1 : k \equiv 0,1,4 \pmod 5 \,\}$ and $S_2 = \{\, k \geq 1 : k \equiv 0,2,3 \pmod 5 \,\}$ for the two sets of admissible parts, and

$$F_{0,1,4}(q) = \sum_{n\geq 0} f_{0,1,4}(n)\, q^n, \qquad F_{0,2,3}(q) = \sum_{n\geq 0} f_{0,2,3}(n)\, q^n$$

for the corresponding generating functions. For a formal series $G(q) = \sum_{n\geq 0} g(n)q^n$ and integers $A \geq 1, B \geq 0$ we write

$$G\,|_{An+B} := \sum_{m\geq 0} g\,(Am+B)\,q^m$$

for the *dissection* of $G$ on the progression $An + B$. (We avoid the more common notation $G(An + B)$, which is easily confused with evaluation.)

Finally, for an odd integer $c$ we set

$$\theta_c(q) \;=\; \sum_{n=-\infty}^{\infty} (-1)^n\, q^{\frac{n(5n+c)}{2}}.$$

**Lemma 2.1.** *Let c be odd. Then every exponent $n(5n + c)/2$ is an integer, and $\theta_{-c} = \theta_c$. If in addition $|c| \leq 5$, then every exponent is non-negative, so $\theta_c \in \mathbb{Z}[[q]]$.*

*Proof.* If $n$ is even then $n(5n + c)$ is even; if $n$ is odd then $5n + c$ is even. Replacing $n$ by $-n$ sends $n(5n + c)/2$ to $n(5n - c)/2$ and fixes $(-1)^n$, which gives $\theta_{-c} = \theta_c$. For the last claim, $5n^2 + cn \geq 0$ for all $n \in \mathbb{Z}$ as soon as $|c| \leq 5$, since the only integers at which $5n^2 + cn$ could be negative are $n = \pm 1$, where its value is $5 \mp |c| \geq 0$. ■

By Lemma 2.1 there are only two series $\theta_c$ with $|c| \leq 5$ that will concern us. We abbreviate

$$T_1(q) := \theta_3(q) = \theta_{-3}(q), \qquad T_2(q) := \theta_1(q) = \theta_{-1}(q).$$

## 3. Reduction modulo 3

**Lemma 3.1.** *In $\mathbb{Z}[[q]]$ we have*

$$F_{0,1,4}(q) \equiv ((q;q^5)_\infty (q^4;q^5)_\infty (q^5;q^5)_\infty)^2 \pmod 3,$$

$$F_{0,2,3}(q) \equiv ((q^2;q^5)_\infty (q^3;q^5)_\infty (q^5;q^5)_\infty)^2 \pmod 3.$$

*Proof.* A part $k$ may occur 0, 1 or 2 times, so

$$F_{0,1,4}(q) = \prod_{k\in S_1} (1 + q^k + q^{2k}) = \prod_{k\in S_1} \frac{1-q^{3k}}{1-q^k}.$$

For a *prime* $\ell$ the binomial theorem gives $(1 - x)^\ell \equiv 1 - x^\ell \pmod{\ell}$, since $\ell$ divides $\binom{\ell}{i}$ for $0 < i < \ell$. Taking $\ell = 3$ and $x = q^k$ yields $1 - q^{3k} \equiv (1 - q^k)^3 \pmod 3$, whence

$$\prod_{k\in S_1} \frac{1-q^{3k}}{1-q^k} \;\equiv\; \prod_{k\in S_1} \frac{(1-q^k)^3}{1-q^k} \;=\; \prod_{k\in S_1} (1-q^k)^2 \pmod 3.$$

Splitting $S_1$ into the residue classes 1, 4 and 0 modulo 5 turns the last product into the asserted product of $q$-Pochhammer symbols. The argument for $F_{0,2,3}$ is identical with $S_2$ in place of $S_1$. ■

**Lemma 3.2 (Jacobi triple product).** *For formal $a$ and $z$,*

$$\sum_{n=-\infty}^{\infty} (-1)^n \, a^{\frac{n(n+1)}{2}} z^n \;=\; (a;a)_\infty \left(\frac{1}{z};a\right)_\infty (za;a)_\infty.$$

This is the Jacobi triple product in the form obtained from the statement in Andrews [1, Ch. 2] by replacing the summation variable $n$ by $-n$.

**Proposition 3.3.** *We have*

$$(q;q^5)_\infty(q^4;q^5)_\infty(q^5;q^5)_\infty = T_1(q), \qquad (q^2;q^5)_\infty(q^3;q^5)_\infty(q^5;q^5)_\infty = T_2(q),$$

*and consequently*

$$F_{0,1,4}(q) \equiv T_1(q)^2, \qquad F_{0,2,3}(q) \equiv T_2(q)^2 \pmod 3.$$

*Proof.* Put $a = q^5$ and $z = q^{-4}$ in Lemma 3.2. Then $(a;a)_\infty = (q^5;q^5)_\infty$, $(1/z;a)_\infty = (q^4;q^5)_\infty$ and $(za;a)_\infty = (q;q^5)_\infty$, while the left-hand side is

$$\sum_{n=-\infty}^{\infty} (-1)^n \, q^{\frac{5n(n+1)}{2}-4n} = \sum_{n=-\infty}^{\infty} (-1)^n \, q^{\frac{n(5n-3)}{2}} = \theta_{-3}(q) = T_1(q).$$

Taking instead $z = q^{-2}$ gives $(1/z;a)_\infty = (q^2;q^5)_\infty$, $(za;a)_\infty = (q^3;q^5)_\infty$ and left-hand side $\theta_1(q) = T_2(q)$. The final assertion now follows from Lemma 3.1. ■

**Remark 3.4.** Multiplying the two identities of Proposition 3.3 and using $(q;q^5)_\infty(q^2;q^5)_\infty(q^3;q^5)_\infty(q^4;q^5)_\infty(q^5;q^5)_\infty = (q;q)_\infty$ gives $T_1T_2 = (q;q)_\infty(q^5;q^5)_\infty$, and hence

$$F_{0,1,4}(q)\, F_{0,2,3}(q) \equiv ((q;q)_\infty(q^5;q^5)_\infty)^2 \pmod 3.$$

So the two functions are, modulo 3, complementary halves of a single eta-quotient. This is one reason to expect them to be exchanged, rather than treated separately, by arithmetic operations on the exponents.

## 4. The dissection lemma

Everything from here on is an analysis of the coefficients of $\theta_c^2$. Expanding the square,

$$\theta_c(q)^2 \;=\; \sum_{j,k=-\infty}^{\infty} (-1)^{j+k} \, q^{E_c(j,k)}, \qquad E_c(j,k) := \frac{j(5j+c)+k(5k+c)}{2}.$$

The whole argument rests on the following exact identity in the integers. Setting

$$X := 10j + c, \qquad Y := 10k + c,$$

we have $X^2 = 20(5j^2 + cj) + c^2$ and likewise for $Y$, so that

$$X^2 + Y^2 \;=\; 40\, E_c(j,k) + 2c^2.$$

No division occurs in (4.2); it is an identity between integers, and all the congruences below are obtained from it by multiplying by units.

**Lemma 4.1.** *Let $p \equiv 3 \pmod 4$ be prime and let $X, Y \in \mathbb{Z}$. Then*

$$p \mid X^2 + Y^2 \Leftrightarrow p \mid X \text{ and } p \mid Y \Leftrightarrow p^2 \mid X^2 + Y^2.$$

*Proof.* Since $p \equiv 3 \pmod 4$, the residue $-1$ is a quadratic non-residue modulo $p$. If $p \mid X^2 + Y^2$ and $p \nmid X$, then $(YX^{-1})^2 \equiv -1 \pmod p$, a contradiction; hence $p \mid X$, and then $p \mid Y^2$, so $p \mid Y$. The remaining implications are immediate. ■

Note that $p \equiv 3 \pmod 4$ forces $p \neq 2$ and $p \neq 5$, so $\gcd(p, 40) = 1$ throughout.

**Lemma 4.2 (the parameters).** *Let $p \equiv 3 \pmod 4$ be prime and let $c \in \{\pm 1, \pm 3\}$.*

(i) *There is a unique $\gamma = \gamma_c(p) \in \{-3,-1,1,3\}$ with $p\gamma \equiv c \pmod{10}$.*

(ii) *$j_c(p) := (p\gamma_c(p) - c)/10$ is an integer, and $j_c(p)$ is the unique residue class modulo $p$ with $p \mid 10j + c$.*

(iii) *$a_c(p) := j_c(p)(5j_c(p) + c)$ is a non-negative integer and satisfies the exact identity*

$$20\, a_c(p) \;=\; \gamma_c(p)^2 p^2 - c^2; \qquad \text{in particular} \quad 20\, a_c(p) \equiv -c^2 \pmod{p^2}.$$

*Proof.* (i) Since $\gcd(p, 10) = 1$, the congruence $p\gamma \equiv c \pmod{10}$ determines $\gamma$ uniquely modulo 10. As $c$ and $p$ are odd, $\gamma$ is odd; as $5 \nmid c$ and $5 \nmid p$, we have $5 \nmid \gamma$. Hence $\gamma \bmod 10 \in \{1,3,7,9\}$, and exactly one representative of that class lies in $\{-3,-1,1,3\}$.

(ii) Integrality is the definition of $\gamma$. If $p \mid 10j + c$ then $10j \equiv -c \equiv 10j_c(p) \pmod p$, and $\gcd(10, p) = 1$ gives $j \equiv j_c(p) \pmod p$.

(iii) Non-negativity is the inequality $5j^2 + cj \geq 0$ of Lemma 2.1, valid since $|c| \leq 5$. For the identity, write $j_0 = j_c(p)$ and $\gamma = \gamma_c(p)$, so that $10j_0 + c = p\gamma$. Then

$$20\, a_c(p) = 20j_0(5j_0 + c) = (10j_0)(10j_0 + 2c) = (p\gamma - c)(p\gamma + c) = \gamma^2 p^2 - c^2. \qquad ■$$

**Theorem 4.3 (Dissection).** *Let $p \equiv 3 \pmod 4$ be prime, let $c \in \{\pm 1, \pm 3\}$, and write $\gamma = \gamma_c(p)$, $j_0 = j_c(p)$, $a = a_c(p)$ as in Lemma 4.2. Then:*

(i) *For all $j, k \in \mathbb{Z}$, the three conditions $E_c(j,k) \equiv a \pmod p$, $p \mid 10j + c$ and $p \mid 10k + c$, and $E_c(j,k) \equiv a \pmod{p^2}$ are equivalent.*

(ii) $\theta_c^2|_{p^2 m + a} \;=\; \theta_\gamma^2$ *as elements of $\mathbb{Z}[[q]]$; that is, for every $m \geq 0$ the coefficient of $q^{p^2 m + a}$ in $\theta_c^2$ equals the coefficient of $q^m$ in $\theta_\gamma^2$.*

(iii) *If* $N \equiv a \pmod{p}$ *but* $N \not\equiv a \pmod{p^2}$*, then the coefficient of* $q^N$ *in* $\theta_c^2$ *is zero, the corresponding index set being empty.*

*Proof.* (i) Put $X = 10j + c$ and $Y = 10k + c$. By (4.2) we have $X^2 + Y^2 = 40E_c(j,k) + 2c^2$, and by Lemma 4.2(iii) we have $40a + 2c^2 = 2\gamma^2p^2$. Subtracting the second from the first,

$$X^2 + Y^2 \;=\; 40(E_c(j,k) - a) + 2\gamma^2p^2.$$

Since $\gcd(40,p) = 1$, (4.3) shows that $p \mid X^2 + Y^2$ if and only if $E_c(j,k) \equiv a \pmod{p}$, and that $p^2 \mid X^2 + Y^2$ if and only if $E_c(j,k) \equiv a \pmod{p^2}$. Lemma 4.1 identifies both of these with the condition $p \mid X$ and $p \mid Y$, which by Lemma 4.2(ii) is the condition $j \equiv k \equiv j_0 \pmod{p}$.

(ii) By (i) the pairs $(j,k)$ contributing to exponents $\equiv a \pmod{p^2}$ are exactly those with $j = pJ + j_0$ and $k = pK + j_0$ for some $J, K \in \mathbb{Z}$, and $(j,k) \mapsto (J,K)$ is a bijection onto $\mathbb{Z}^2$. For such a pair,

$$X = 10(pJ + j_0) + c = 10pJ + p\gamma = p\,(10J + \gamma), \qquad Y = p\,(10K + \gamma),$$

so by (4.2) applied twice, once with $c$ and once with $\gamma$,

$$40E_c(j,k) + 2c^2 = X^2 + Y^2 = p^2[(10J+\gamma)^2 + (10K+\gamma)^2] = p^2\big[40E_\gamma(J,K) + 2\gamma^2\big].$$

Using $20a = \gamma^2p^2 - c^2$ from Lemma 4.2(iii), the constant terms give $2c^2 + 40a = 2\gamma^2p^2$, and therefore

$$E_c(j,k) \;=\; p^2\,E_\gamma(J,K) + a.$$

Finally, since $p$ is odd,

$$(-1)^{j+k} = (-1)^{pJ+pK+2j_0} = (-1)^{J+K}.$$

Comparing (4.1), (4.4) and (4.5), the coefficient of $q^{p^2m+a}$ in $\theta_c^2$ is $\sum(-1)^{J+K}$ over all $(J,K)$ with $E_\gamma(J,K) = m$, which is the coefficient of $q^m$ in $\theta_\gamma^2$. (That $E_\gamma(J,K) \geq 0$, so that no exponent is lost, is Lemma 2.1 together with $|\gamma| \leq 3$.)

(iii) Immediate from (i): a pair $(j,k)$ with $E_c(j,k) \equiv a \pmod{p}$ automatically has $E_c(j,k) \equiv a \pmod{p^2}$, so if $N \equiv a \pmod{p}$ and $N \not\equiv a \pmod{p^2}$ no pair contributes and the coefficient is an empty sum. ■

## 5. The main theorem

Only the two values $c = -3$ and $c = 1$ are needed, since $F_{0,1,4} \equiv \theta_{-3}^2$ and $F_{0,2,3} \equiv \theta_1^2$ modulo 3 by Proposition 3.3. We therefore abbreviate

$$a(p) := a_{-3}(p), \qquad b(p) := a_1(p),$$

so that, by Lemma 4.2(iii),

$$20\,a(p) = \gamma_{-3}(p)^2p^2 - 9, \qquad 20\,b(p) = \gamma_1(p)^2p^2 - 1.$$

**Theorem 5.1.** *Let $p \equiv 3 \pmod 4$ be prime. Then for all $m \geq 0$:*

*(i) if $p \equiv \pm 1 \pmod 5$,*

$$f_{0,1,4}\big(p^2 m + a(p)\big) \equiv f_{0,1,4}(m), \qquad f_{0,2,3}\big(p^2 m + b(p)\big) \equiv f_{0,2,3}(m) \pmod 3;$$

*(ii) if $p \equiv \pm 2 \pmod 5$,*

$$f_{0,1,4}\big(p^2 m + a(p)\big) \equiv f_{0,2,3}(m), \qquad f_{0,2,3}\big(p^2 m + b(p)\big) \equiv f_{0,1,4}(m) \pmod 3.$$

*Moreover, by quadratic reciprocity, case (i) occurs precisely when 5 is a quadratic residue modulo $p$.*

*Proof.* By Proposition 3.3 and Theorem 4.3(ii),

$$F_{0,1,4}|_{p^2m+a(p)} \equiv \theta_{-3}^2|_{p^2m+a(p)} = \theta_{\gamma_{-3}(p)}^2, \qquad F_{0,2,3}|_{p^2m+b(p)} \equiv \theta_1^2|_{p^2m+b(p)} = \theta_{\gamma_1(p)}^2 \pmod 3.$$

By Lemma 2.1, $\theta_\gamma^2 = T_1^2$ when $|\gamma| = 3$ and $\theta_\gamma^2 = T_2^2$ when $|\gamma| = 1$, and by Proposition 3.3 these are congruent modulo 3 to $F_{0,1,4}$ and $F_{0,2,3}$ respectively. It therefore remains only to determine $|\gamma_{-3}(p)|$ and $|\gamma_1(p)|$.

By Lemma 4.2(i), $\gamma_{-3}(p) \equiv -3\,p^{-1} \pmod{10}$. Since 3 is invertible modulo 10,

$$|\gamma_{-3}(p)| = 3 \Leftrightarrow -3p^{-1} \equiv \pm 3 \pmod{10} \Leftrightarrow p^{-1} \equiv \mp 1 \pmod{10} \Leftrightarrow p \equiv \mp 1 \pmod{10},$$

and for odd $p$ the condition $p \equiv \pm 1 \pmod{10}$ is the same as $p \equiv \pm 1 \pmod 5$. Likewise $\gamma_1(p) \equiv p^{-1} \pmod{10}$, so $|\gamma_1(p)| = 1$ if and only if $p \equiv \pm 1 \pmod 5$. Thus in case (i) we have $|\gamma_{-3}| = 3$ and $|\gamma_1| = 1$, giving $T_1^2$ and $T_2^2$; in case (ii) we have $|\gamma_{-3}| = 1$ and $|\gamma_1| = 3$, giving $T_2^2$ and $T_1^2$, that is, the two functions are exchanged.

For the last sentence, $p \equiv \pm 1 \pmod 5$ says exactly that $p$ is a quadratic residue modulo 5, i.e. $\left(\frac{p}{5}\right) = 1$; since $5 \equiv 1 \pmod 4$, quadratic reciprocity gives $\left(\frac{p}{5}\right) = \left(\frac{5}{p}\right)$. ■

**Corollary 5.2 (closed forms for the constants).** *Let $p \equiv 3 \pmod 4$ be prime. If $p \equiv \pm 1 \pmod 5$ then*

$$a(p) = \frac{9\,(p^2-1)}{20}, \qquad b(p) = \frac{p^2-1}{20},$$

*while if $p \equiv \pm 2 \pmod 5$ then*

$$a(p) = \frac{p^2-9}{20}, \qquad b(p) = \frac{9p^2-1}{20}.$$

*Proof.* Immediate from (5.1) and the values of $|\gamma_{-3}(p)|$, $|\gamma_1(p)|$ computed in the proof of Theorem 5.1. ■

The table below lists the data for the primes $p \equiv 3 \pmod 4$ with $p < 80$. Here "exchanged" and "preserved" refer to the two cases of Theorem 5.1, and $d$ is the exponent appearing in Theorem 7.1.

| $p$ | $p$ mod 5 | $\gamma_{-3}$ | $\gamma_1$ | $a(p)$ | $b(p)$ | behavior | $d$ |
|---|---|---|---|---|---|---|---|
| 3 | 3 | −1 | −3 | 0 | 4 | exchanged | 4 |
| 7 | 2 | 1 | 3 | 2 | 22 | exchanged | 4 |
| 11 | 1 | −3 | 1 | 54 | 6 | preserved | 2 |
| 19 | 4 | 3 | −1 | 162 | 18 | preserved | 2 |
| 23 | 3 | −1 | −3 | 26 | 238 | exchanged | 4 |
| 31 | 1 | −3 | 1 | 432 | 48 | preserved | 2 |
| 43 | 3 | −1 | −3 | 92 | 832 | exchanged | 4 |
| 47 | 2 | 1 | 3 | 110 | 994 | exchanged | 4 |
| 59 | 4 | 3 | −1 | 1566 | 174 | preserved | 2 |
| 67 | 2 | 1 | 3 | 224 | 2020 | exchanged | 4 |
| 71 | 1 | −3 | 1 | 2268 | 252 | preserved | 2 |

**Corollary 5.3 (the two smallest cases).** *For all $m \geq 0$,*

$$f_{0,1,4}(9m) \equiv f_{0,2,3}(m), \qquad f_{0,2,3}(9m+4) \equiv f_{0,1,4}(m) \pmod 3,$$

$$f_{0,1,4}(49m+2) \equiv f_{0,2,3}(m), \qquad f_{0,2,3}(49m+22) \equiv f_{0,1,4}(m) \pmod 3.$$

*Proof.* Take $p = 3$ and $p = 7$ in Theorem 5.1(ii). For $p = 3$: $\gamma_{-3}(3) = -1$ since $3 \cdot (-1) \equiv -3 \pmod{10}$, so $j_{-3}(3) = 0$ and $a(3) = 0$; and $\gamma_1(3) = -3$ since $3 \cdot (-3) \equiv 1 \pmod{10}$, so $j_1(3) = -1$ and $b(3) = (-1)(-5+1) = 4$. For $p = 7$: $\gamma_{-3}(7) = 1$, $j_{-3}(7) = 1$, $a(7) = 1 \cdot 2 = 2$; and $\gamma_1(7) = 3$, $j_1(7) = 2$, $b(7) = 2 \cdot 11 = 22$. ■

It may be worth unwinding the case $p = 3$ of Theorem 4.3 explicitly, since it is the shortest instance of the whole argument. With $c = -3$ we have $3 \mid 10j - 3$ if and only if $3 \mid j$; writing $j = 3J$ and $k = 3K$ gives $j(5j-3) = 9J(5J-1)$, hence $E_{-3}(j,k) = 9\,E_{-1}(J,K)$, while $(-1)^{j+k} = (-1)^{J+K}$. Thus $\theta_{-3}^2|_{9m} = \theta_{-1}^2 = T_2^2$, which is $F_{0,2,3}$ modulo 3.

## 6. Congruences with exceptions

Theorem 4.3(iii) says that the surviving terms of $\theta_c^2$ in a progression of modulus $p$ are already confined to a single progression of modulus $p^2$. This yields congruences with exceptions in the sense of Keith [5]: the mechanism is not cancellation but vacuity, the coefficient outside the exceptional class being an empty sum.

**Theorem 6.1.** *Let* $p \equiv 3 \pmod 4$ *be prime and write*

$$a(p) = p\alpha + a_0, \quad 0 \le a_0 < p, \qquad b(p) = p\beta + b_0, \quad 0 \le b_0 < p.$$

*Then for all* $m \ge 0$*,*

$$f_{0,1,4}(pm + a_0) \equiv 0 \pmod 3 \quad \text{whenever } m \not\equiv \alpha \pmod p,$$

$$f_{0,2,3}(pm + b_0) \equiv 0 \pmod 3 \quad \text{whenever } m \not\equiv \beta \pmod p.$$

*Proof.* Put $N = pm + a_0$, so $N \equiv a(p) \pmod p$. Then $N \equiv a(p) \pmod{p^2}$ if and only if $pm \equiv p\alpha \pmod{p^2}$, that is, if and only if $m \equiv \alpha \pmod p$. If $m \not\equiv \alpha \pmod p$, then by Theorem 4.3(iii) the coefficient of $q^N$ in $\theta_{-3}^2$ vanishes, and by Proposition 3.3 this coefficient is congruent to $f_{0,1,4}(N)$ modulo 3. The second statement is the same argument with $c = 1$. ■

**Corollary 6.2.** *For all* $m \ge 0$*,*

$$f_{0,1,4}(3m) \equiv 0 \pmod 3 \text{ unless } 3 \mid m, \qquad f_{0,2,3}(3m+1) \equiv 0 \pmod 3 \text{ unless } m \equiv 1 \pmod 3,$$

$$f_{0,1,4}(7m+2) \equiv 0 \pmod 3 \text{ unless } 7 \mid m, \qquad f_{0,2,3}(7m+1) \equiv 0 \pmod 3 \text{ unless } m \equiv 3 \pmod 7.$$

*Proof.* Apply Theorem 6.1 with $a(3) = 0$, $b(3) = 4 = 3 \cdot 1 + 1$, $a(7) = 2$ and $b(7) = 22 = 7 \cdot 3 + 1$. ■

These congruences are sharp: the exceptional classes are exactly the progressions on which Theorem 5.1 places the surviving terms, and on those progressions the values are in general not divisible by 3. For instance, among the 1225 values $f_{0,1,4}(49t + 2)$ with $t \le 1224$, exactly 476 are not divisible by 3.

## 7. Self-similarity

**Theorem 7.1.** *Let* $p \equiv 3 \pmod 4$ *be prime and set*

$$d = d(p) = \begin{cases} 2, & p \equiv \pm 1 \pmod 5, \\ 4, & p \equiv \pm 2 \pmod 5. \end{cases}$$

*Then* $p^d \equiv 1 \pmod{20}$*, and for all* $n \ge 0$ *and* $m \ge 0$*,*

$$f_{0,1,4}\left(p^{dn} m + \frac{9\,(p^{dn} - 1)}{20}\right) \equiv f_{0,1,4}(m), \qquad f_{0,2,3}\left(p^{dn} m + \frac{p^{dn} - 1}{20}\right) \equiv f_{0,2,3}(m) \pmod 3.$$

*In particular both constants are integers.*

*Proof.* Since $p$ is odd, $p^2 \equiv 1 \pmod 8$, so $p^2 \equiv 1 \pmod 4$. Modulo 5: if $p \equiv \pm 1$ then $p^2 \equiv 1$; if $p \equiv \pm 2$ then $p^2 \equiv 4$ and $p^4 \equiv 1$. In either case $p^d \equiv 1 \pmod 4$ and $p^d \equiv 1 \pmod 5$, hence $p^d \equiv 1 \pmod{20}$ and $20 \mid p^{dn} - 1$ for all $n$.

Write $A_n = 9(p^{dn} - 1)/20$ and $B_n = (p^{dn} - 1)/20$. These satisfy $A_0 = B_0 = 0$ and

$$A_n = p^d A_{n-1} + A_1, \qquad B_n = p^d B_{n-1} + B_1 \qquad (n \geq 1).$$

We first check the case $n = 1$.

If $p \equiv \pm 1 \pmod 5$ then $d = 2$ and, by Corollary 5.2, $A_1 = 9(p^2 - 1)/20 = a(p)$ and $B_1 = (p^2 - 1)/20 = b(p)$, so the case $n = 1$ is exactly Theorem 5.1(i).

If $p \equiv \pm 2 \pmod 5$ then $d = 4$, and by Corollary 5.2 we have $20a(p) = p^2 - 9$ and $20b(p) = 9p^2 - 1$, whence

$$20\big(p^2 b(p) + a(p)\big) = p^2(9p^2 - 1) + (p^2 - 9) = 9(p^4 - 1), \qquad 20\big(p^2 a(p) + b(p)\big) = p^2(p^2 - 9) + (9p^2 - 1) = p^4 - 1.$$

Thus $A_1 = p^2 b(p) + a(p)$ and $B_1 = p^2 a(p) + b(p)$. Applying Theorem 5.1(ii) twice,

$$f_{0,1,4}(p^4 m + A_1) = f_{0,1,4}\left(p^2\big(p^2 m + b(p)\big) + a(p)\right) \equiv f_{0,2,3}\big(p^2 m + b(p)\big) \equiv f_{0,1,4}(m) \pmod 3,$$

and symmetrically for $f_{0,2,3}$ with $B_1$.

The general case follows by induction on $n$: by (7.1),

$$p^{dn} m + A_n = p^d\big(p^{d(n-1)} m + A_{n-1}\big) + A_1,$$

so the case $n = 1$ gives $f_{0,1,4}(p^{dn} m + A_n) \equiv f_{0,1,4}\big(p^{d(n-1)} m + A_{n-1}\big)$, which is $\equiv f_{0,1,4}(m)$ by the induction hypothesis. The case $n = 0$ is trivial. ■

For $p = 3$ this gives the family on progressions of modulus $81^n$,

$$f_{0,1,4}\left(81^n m + \frac{9(81^n - 1)}{20}\right) \equiv f_{0,1,4}(m), \qquad f_{0,2,3}\left(81^n m + \frac{81^n - 1}{20}\right) \equiv f_{0,2,3}(m) \pmod 3,$$

whose first non-trivial instances are $f_{0,1,4}(81m + 36) \equiv f_{0,1,4}(m)$ and $f_{0,2,3}(81m + 4) \equiv f_{0,2,3}(m)$. For $p = 7$ the modulus is $2401^n$ and the constants are $9(2401^n - 1)/20$ and $(2401^n - 1)/20$, the first instances being $f_{0,1,4}(2401m + 1080) \equiv f_{0,1,4}(m)$ and $f_{0,2,3}(2401m + 120) \equiv f_{0,2,3}(m)$. For $p = 11$, where the two functions are preserved rather than exchanged, the modulus is already $121^n$: $f_{0,1,4}(121m + 54) \equiv f_{0,1,4}(m)$ and $f_{0,2,3}(121m + 6) \equiv f_{0,2,3}(m)$.

## 8. Concluding remarks

**8.1. The modulus 3 comes from the multiplicity bound, and the method does not survive raising it.** The prime 3 in every congruence above originates not in the residue condition on the parts but in the bound on multiplicity, through the factorisation $1 + q^k + q^{2k} = (1 - q^{3k})/(1 - q^k)$. If each part may be used at most $t$ times and $\ell := t + 1$ is prime, the same reduction gives

$$\prod_{k \in S} \frac{1 - q^{\ell k}}{1 - q^k} \equiv \prod_{k \in S} (1 - q^k)^{\ell - 1} \pmod{\ell},$$

so the generating function becomes an $(\ell - 1)$-st power of a theta series. For $\ell = 3$ this is a square, and the dissection reduces to the classical fact that a sum of *two* squares vanishing modulo a prime $p \equiv 3 \pmod 4$ must vanish termwise. For $\ell = 5$ one obtains a fourth power, and the corresponding statement about sums of four squares is false for every prime, since every residue class is a sum of four squares in many ways; the argument of Section 4 therefore has no analogue there, and any congruences for $t = 4$ would have to be found by other means. For $\ell = 2$ (distinct parts) the reduction gives a single theta series modulo 2, where a one-square argument is available. If $\ell$ is composite, say $\ell = \ell' s$ with $s > 1$, the reduction is not even clean: one obtains $\prod (1 - q^{sk})^{\ell'}/(1 - q^k)$ modulo $\ell'$, which does not collapse to a power of a single theta series. The restriction to $t = 2$ is thus not an accident of exposition.

**8.2. The role of $p \equiv 3 \pmod 4$.** Lemma 4.1 is the only place where the hypothesis on $p$ is used, and it is essential: for $p \equiv 1 \pmod 4$ the residue $-1$ is a square, so a sum of two squares can vanish modulo $p$ with neither term vanishing, and the set of pairs $(j, k)$ contributing to a given progression is neither empty nor of the required shape. A numerical search over all residues modulo $p^2$ for $p = 5,13,17,29,37$ found no congruence of the form of Theorem 5.1 in any of the four possible pairings, which suggests the obstruction is real and not merely an artefact of the method.

**8.3. Exchange versus preservation.** The dichotomy in Theorem 5.1 is governed by $\left(\frac{5}{p}\right)$, and its source is visible in Lemma 4.2(i): the parameter $\gamma_c(p) \equiv c\, p^{-1} \pmod{10}$ records how multiplication by $p$ permutes the classes $\pm 1$ and $\pm 2$ modulo 5, which are precisely the two Rogers–Ramanujan classes indexing $T_2$ and $T_1$. The exchange therefore happens exactly when $p$ is a non-residue modulo 5. It also explains the difference between the exponents $d = 2$ and $d = 4$ in Theorem 7.1: in the exchanging case one must apply the dissection twice to return to the same function, which is why the self-similarity of $f_{0,1,4}$ appears at $81^n$ and $2401^n$ rather than at $9^n$ and $49^n$.

**8.4. Questions.** Three questions seem natural. First, the congruences do not lift: a direct computation shows that $f_{0,1,4}(9m) \equiv f_{0,2,3}(m)$ and $f_{0,1,4}(49m + 2) \equiv f_{0,2,3}(m)$ both fail modulo 9, so the phenomenon is genuinely a modulo-3 one, and it would be interesting to know what the correct modulo-9 statement is, if any. Second, is there an analogous pair of functions attached to a prime other than 5 — that is, a partition of the non-zero residues modulo some prime $r$ into two classes, stable under $x \mapsto -x$, whose restricted partition functions with multiplicity at most two satisfy congruences of this type? The proof above uses $r = 5$ only through the exact identity (4.2); the analogue of that identity for a general $r$ has the form $X^2 + Y^2 = 8rE + 2c^2$ with $X = 2rj + c$, and it is the requirement $|c| \leq r$ of Lemma 2.1 that limits the possible $\gamma$. Third, one may ask for a combinatorial or bijective explanation of the exchange in Theorem 5.1(ii); the proof given here is entirely generating-function-theoretic and offers none.

## 9. Acknowledgements

The author initiated and directed this investigation, curated the intermediate results, chose which statements to develop, and wrote and edited the final exposition. The original pair of partition functions was suggested by an AI language model (Anthropic Claude), which was also used to help draft parts of the argument, to carry out the numerical verifications reported above, and to generalize an earlier version of Theorem 5.1 that treated only the case $p = 7$. Every definition, statement and proof in the present paper has been checked by hand by the author.